\documentclass[paper]{amsart}
\usepackage{amssymb,amscd}
\usepackage{fullpage}

\newtheorem{thm}{Theorem}[section]

\newtheorem{lem}[thm]{{Lemma}}
\newtheorem{defn}[thm]{Definition}

\numberwithin{equation}{section}

\def\overlap{\operatorname{overlap}}

\def\Z{\mathbb Z}

\def\Mod{\operatorname{Mod}}

\def\comp{\operatorname{comp}}
\def\smg{\operatorname{smg}}

\title[GROEBNER-SHIRSHOV BASES]
{NOTE ON THE CALCULATION OF GROEBNER-SHIRSHOV BASES FOR  AFFINE
WEYL GROUPS}
\author{Cenap \"{O}zel$^{1}$, Adem K{\i}l{\i}\c cman$^{2}$ \& Erol Yilmaz$^{1}$}
\address{$^{1,3}$ AIBU Golkoy Kampusu, Bolu 14280, Turkey.}
\email{$^{1,3}$ cenap@ibu.edu.tr\  \, ; yerol@ibu.edu.tr}
\address{Department of Mathematics and Institute for Mathematical Research, University Putra Malaysia, 43400 UPM,
Serdang, Selangor, Malaysia} \email{$^{2}$
akilicman@putra.upm.edu.my}
\begin{document}

\begin{abstract}
\noindent In this work we will consider the calculation of
Groebner-Shirshov bases of Coxeter groups. This will be the main
focus of the work. In \cite{Bokut-Shiao}, Bokut \& Shiao gave the
Groebner-Shirshov bases of positive definite classical Coxeter
groups $A_l, B_l, D_l$ by using the techniques of Elimination of
Leading Word.

\noindent We will give a counter example to a hypothesis which is
introduced by Bokut \& Shiao in \cite{Bokut-Shiao} and we will
calculate the Groebner-Shirshov bases of the positive degenerate
infinite affine Weyl group $\widetilde{A}_n $ which is isomorphic
to semi-direct product group $\Sigma_n \ltimes {\mathbb Z}^{n-1}$,
and further we classify all the reduced elements of the group by
using the Composition Diamond Lemma.
\end{abstract}

\maketitle

\section{Introduction}

\noindent The Grobner basis theory for commutative algebras was
introduced by Buchberger \cite{Buchberger} and provides a solution
to the reduction problem for commutative algebras.  It is also
known an effective algorithm of computing a set of generators for
a given ideal of a commutative ring which further can be used to
determine the reduced elements with respect to the relations given
by the ideal. Later in \cite{Bergman}, Bergman generalized the
Grobner basis theory to associative algebras by proving the {\it
Diamond Lemma} which is the key ingredient in the theory and also
equivalent to the so-called {\it Composition Lemma} which
characterizes the leading terms of elements in the given ideal,
see \cite{ufnarovski}.\\

\noindent First of all we recall some facts about non-commutative
Grobner bases which are also known in the literature as
Grobner–-Shirshov bases (see, for example, \cite{Newman} and
\cite{Shirshov}).\\

\noindent Let $S$ be a linearly ordered set, $k$ is a field,
$k\langle S\rangle$ is the formal free associative algebra over
$S$ and $k$ where generators are the letters in the totally
ordered set $S$. Further we let $S^{*}$ be the set of words which
are obtained from letters in $S$ performing free product.\\

\noindent Now suppose that the words in $S^{*}$ are
degree-lexicographic ordered.i.e. comparing two words first by
lengths and then lexicographically.
\begin{defn}
Any polynomial $f\in k\langle S\rangle$ is a linear sum of words in
$S^{*}$.
\end{defn}
In this section the words plays the crucial roles similar to the
role of monomials. Since any polynomial in $f\in k\langle
S\rangle$ consists of finitely many words and the
degree-lexicographic order is total, any polynomial $f\in k\langle
S\rangle$ has a leading word. Let $\bar{f}$ be the leading word of
$f$.

\begin{defn}
If $\bar{f}$ in $f$ has a coefficient $1$ then we say that $f$ is
called  \emph{monic}.
\end{defn}
\begin{defn}[\emph{Composition of Intersection}]
A composition of intersection  of two monic polynomials relative to
some word $w$ is defined as $(f,g)_{w} = fb-ag$ whenever
$w=\bar{f}b=a\bar{g}$ , $deg(\bar{f}) + deg(\bar{g})> deg(w)$.
\end{defn}

\begin{defn}[\emph{Composition of Including}]
A composition of including $(f,g)_{w}$ of two monic polynomials
relative to some word $w$ is defined as $(f,g)_{w} = f-agb$ whenever
$w=\bar{f}=a\bar{g}b$.
\end{defn}

\begin{defn}[\emph{Elimination of the Leading Word }]
In the last case $(f,g)_{w}=f-agb $ is called the \emph{elimination
of the leading word (ELW) of $g$ in $f$}.
\end{defn}

\begin{defn}[\emph{Trivial Relative}]
A composition $(f,g)_{w}$ is called \emph{trivial relative} to some
$R\subset k \langle X\rangle$ if $(f,g)_{w}= \sum
\alpha_{i}a_{i}t_{i}b_{i}$, where $t_{i}\in R$, $a_{i},b_{i}\in
S^{*}$ and $a_{i}{\bar{t}}_{i}b_{i}< w $.
\end{defn}
\noindent In particular, if $(f,g)_{w}$ goes to zero by the $ELW$
's of $R$ then $(f,g)_{w}$ is trivial relative to $R$.

\begin{defn}[\emph{Groebner-Shirshov Bases}]
A \emph{Groebner-Shirshov basis} is a subset $R$ of $k\langle S
\rangle$ if any composition of polynomials from $R$ is trivial
relative to $R$.
\end{defn}

\noindent By the algebra $\langle S|R\rangle$ with generators $S$
and defining relations $R$, we mean the \emph{factor algebra} of
$k\langle S \rangle$ by the ideal generated by $R$. The following
lemma goes back to the \emph{Diamond Lemma} of \cite{Newman} and
is also known as the \emph{Composition Lemma} of \cite{Shirshov}.

\begin{lem}[\emph{Composition-Diamond Lemma}]
$R$ is a Groebner-Shirshov basis if and only if the set
$$\{u\in S^{*} | u\neq a\bar{f}b,\; \forall f\in
R\}$$ of $R$ reduced words, consists of a linear basis of the
algebra $\langle S | R \rangle$.
\end{lem}

\noindent Note that if a subset $R$ of $k\langle S \rangle$ is not
a Groebner-Shirshov basis then we can add all non-trivial
compositions of polynomials of $R$ to $R$  and by continuing this
process (possibly infinitely) many times in order to get a
Groebner-Shirshov basis $S^{\comp}$ that contains $R$. This
procedure is called \emph{Buchberger-Shirshov Algorithm}, for more
details, see \cite{Buchberger}.

\begin{defn}[\emph{Reduced Groebner-Shirshov Basis}]
A Groebner-Shirshov basis is called \emph{reduced }if any $s\in S$
is a linear combination of $S \setminus \{ \{s\} - \mbox{ reduced
words }\}$.
\end{defn}
We also note that any ideal of $k\langle S \rangle$ has a unique
reduced Groebner-Shirshov basis. Now let $R$ be a set of
\emph{semigroup relations} i.e.
$$
\{u-v : \, u,v\in S^{*}\},
$$
then any nontrivial composition will have the same form. Thus the
set $s^{\comp}$ also consists of semigroup relations.\\

\noindent Let $A= \smg \langle S|R\rangle$ be a semigroup
presentation. Then $R$ is a subset of $k\langle S \rangle$ and we
can find a Groebner-Shirshov basis $s^{\comp}$. The last set doesn't
depend on $k$, and it consist of semigroup relations. we will call
$R^{comp}$ a Groebner-Shirshov basis of $A$. It is the same as
Groebner-Shirshov basis of the semigroup algebra $kA= \langle S |
R\rangle$.\\

\noindent The same terminology is valid for any group presentation
meaning that we include in this presentation all trivial group
relations of the form $ss^{-1}= 1$, $s^{-1}s= 1$, $s\in S$.\\

\noindent In the following sections we will consider the calculation
of Groebner-Shirshov basis of Weyl groups and this will be the main
focus of the present work. In \cite{Bokut-Shiao}, Bokut \& Shiao
gave the Groebner-Shirshov basis of positive definite classical Weyl
groups $A_l, B_l, D_l$ by using the techniques of \emph{Composition
Diamond Lemma} and \emph{Elimination of Leading Word}. In the next
section we will give the Groebner-Shirshov bases of Weyl groups
$A_l, B_l, D_l$. In the last two sections of the work we will give a
counter-example to the hypothesis that was introduced by Bokut \&
Shiao in \cite{Bokut-Shiao} and later we will calculate the
Groebner-Shirshov basis of the positive degenerate infinite affine
Weyl group $\widetilde{A}_n $ which is isomorphic to semi-direct
product group $\Sigma_n \ltimes \Z^{n-1}$. Finally, we give reduced
elements of $\widetilde{A}_n $.

\section{Calculations for finite Weyl Groups}
From \cite{Bokut-Shiao} we have the following results.
\subsection{The Weyl Group $A_{l}$.}
\begin{thm}\label{cenap}
The Weyl group $A_{l}$ is generated by $\{s_{1},s_{2},\cdots,
s_{l}\}$ and $s_{i,j} = s_{i}s_{i-1}\ldots s_{j}$, where $i>j$;
and $s_{i,i} = s_{i}$ , $s_{i,i+1} = 1$ with defining relations;

\begin{tabular}{llll}
$A{1}$ :& ${s_{i}}^2= 1$ $ \dashrightarrow $ $s_{i}s_{i}=1$ & where
& $1
\leq i < l$, \\ $A{2}$ : & $s_{i}s_{j} = s_{j}s_{i}$ & where & $i-j>1$, \\
$A{3}$ : & $s_{i+1}s_{i}s_{i+1} = s_{i}s_{i+1}s_{i}$ & where &
$1\leq i\leq l-1$, \\ $A{4}$ : & $s_{i+1,j}s_{i+1} = s_{i}s_{i+1,j}$
&
where & $1 \leq j \leq i$.\\
\end{tabular}

\noindent So the Groebner-Shirshov basis for $A_{l}$ consists of the
relations that which was given by $A{1}-A{4}$.
\end{thm}

\subsection{The Weyl Group $ B_{l}$.}
\begin{thm} The Weyl group $B_l$ is generated by $s_{i}$, $1 \leq i \leq l$ and
$s_{i,j} = s_{i}s_{i-1}\ldots s_{j},$ where $i>j$; and $s_{i,i} =
s_{i}$ , $s_{i,i+1} = 1$ with defining relations;
\begin{tabular}{llll}
$B{1}$: & ${s_i}^{2}= 1$ $ \dashrightarrow  s_{i}s_{i}=1$, & where &
$1 \leq i \leq l$ and $i-j > 1$, \\
$B{2}$: & $s_{i}s_{j} = s_{j}s_{i}$, & where & $i-j>1$, \\ $B{3}$:&
$s_{i+1}s_{i}s_{i+1} = s_{i}s_{i+1}s_{i}$, & where & $1 \leq i \leq
l-2$, \\ $B{4}$: & $s_{i+1,j}s_{i+1} = s_{i}s_{i+1,j}$, & where & $j
\leq i \leq l-2$, \\
$B{5}$: & $s_{l}s_{l-1}s_{l}s_{l-1} = s_{l-1}s_{l}s_{l-1}s_{l}$ & \\
$B{6}$:& $s_{l,j}s_{l,j} = s_{l-1}s_{l,j}s_{l,j+1}$, & where
& $j<l$.\\
\end{tabular}

\noindent So the Groebner-Shirshov basis for $B_{l}$ consists of the
relations that given by $B{1}-B{6}$.
\end{thm}

\subsection{The Weyl Grup $ D_{l}$.}
\begin{thm}
The Weyl group $ D_{l}$ is generated by $s_{i}$ , $1 \leq i \leq l$
together with defining relations $ A_{l-1}$ and $s_{ij}$ , $1 \leq j
\leq i+1 \leq l$; and $s_{l,j} = s_{l}s_{l-2}\ldots s_{j},$ where $j
\leq l-2$; and $s_{l,l} = 1$ , $s_{l,l-1} = s_{l}$.

$D{1}$: $s_{l,l}= 1$ $ \dashrightarrow $ $s_{l}s_{l}=1$,

$D{2}$: $s_{l}s_{l-1} = s_{l-1}s_{l}$,

$D{3}$: $s_{l}s_{l-2}s_{l} = s_{l-2}s_{l}s_{l-2}$,

$D{4}$ : $s_{l,j}s_{l-1,j} = s_{l-1}s_{l,j}s_{l-1,j+1}$ where $j
\leq l-2$,

$D{5}$: $s_{l,j}s_{l-1}s_{l} = s_{l-2}s_{l,j}s_{l-1}$ where $j \leq
l-2$,

$D{6}$: $s_{l,j}s_{l-1,k} = s_{l-2}s_{l,j}s_{l-1,k}s_{l,k+1}$ where
$j<k \leq l-2$,

$D{7}$: $s_{i,i}= 1$ $ \dashrightarrow s_{i}s_{i}=1$ where $1 \leq i
\leq l-1$,

$D{8}$: $s_{i}s_{j} = s_{j}s_{i}$ where $i-j>2$,

$D{9}$: $s_{i+1}s_{i}s_{i+1} = s_{i}s_{i+1}s_{i}$ where $1\leq i\leq
l-2$,

$D{10}$ : $s_{i+1,j}s_{i+1} = s_{i}s_{i+1,j}$.

\noindent So the Groebner-Shirshov basis for $D_{l}$ consists of the
relations given by $D{1}-D{10}$ and relations of $A_{l-1}$.
\end{thm}

\section{A Counter-example to Hypothesis of Bokut\& Shiao.}

\noindent In \cite{Bokut-Shiao} they formulate a general hypothesis
on Groebner-Shirshov bases for any Coxeter group. In the next we
will provide a counter example.\\

\noindent Let $W$ be a Coxeter group, $S$ is a linear ordered set
with $l$-generators, and $M=[m_{ss'}]$ be a $l\times l$ Coxeter
matrix associated with $W$. Then we have
$$
W=\smg\langle s^2=1,\;(ss')^{m_{ss'}}=1,\; s\neq s',\; \forall s,s'
\in S \; \mbox{and finite}\;m_{ss'}\rangle. $$  Now define for
finite $m_{ss'}$,
\begin{align}
m(s,s')&=ss' ...\quad(\mbox{there are $m_{ss'}$ alternative letters
$s,s'$}),\notag\\
(m-1)(s,s')&=ss' ...\quad(\mbox{there are $m_{ss'}-1$ alternative
letters $s,s'$}),\notag
\end{align}
and so on. In particular, if $m_{ss'}=2$ then
$m(s,s')=ss',(m-1)(s,s')=s$. Similarly, if $m_{ss'}=3$, then
$m(s,s')=ss's, (m-1)(s,s')=ss'$ where $m_{ii} = 1$ that is $(s_i)^2
= 1$ for all $i$; thus the generators are involutions.\\

Now if $m_{ij} = 2$, then the generators $s_i$ and $s_j$ commute
thus the notation to define the relations of $W$ can be presented in
the form
$$s^2=1,\,m(s,s')=m(s',s),\; s>s',\quad(*)
$$
for all $s,s' \in S$ and finite $m_{ss'}$.

\begin{defn}
Two words in $S$ are \emph{equivalent} if they are equal module
commutativity relations from $(*)$.
\end{defn}
To be more precise, it means that they are equal in the so called
free partially commutative semigroup (algebra)generated by $S$ and
with commutativity relations from $(*)$. Later on the word problem
can be solved in this semigroup(or algebra).

\begin{defn}
Two relations $a=b$ and $c=d$ of $W$ are called \emph{equivalent} if
$a,b$ are equivalent to $c,d$ respectively.
\end{defn}

\subsection{The Hypothesis.} The Groebner-Shirshov bases of $W$ consist of
initial relations $(*)$ and relations that are equivalent to the
following ones:

$$
(m-1)(s,s')(m-1)(s_{i_{1}},s_{i_{2}})...(m-1)(s_{i_{2k-1}},s_{i_{2k}})
m(s_{i_{2k+1}}s_{i_{2k+2}})
$$
$$
=m(s',s)(m-1)(s_{i_{1}},s_{i_{2}})...(m-1)(s_{i_{2k-1}},s_{i_{2k}})(m-1)(s_{i_{2k+1}}
s_{i_{2k+2}}),\quad (**)
$$
where $s>s',\, s_{i_{1}} < s_{i_{2}},...,s_{i_{2k-1}}
<s_{i{2k}},s_{i_{2k+1}} < s_{i_{2k+2}}$, and any neighbor pairs
$(s',s)(s_{i_{1}},s_{i_{2}}),...,(s_{i_{2k+1}}s_{i_{2k+2}})$ are
different, and
$$
\begin{cases}
s_{i_{2}} =s', &\quad \mbox{if}\; m_{ss'}\; \mbox{is even}\\
s_{i_{2}} =s , &\quad \mbox{if}\; m_{ss'}\; \mbox{is odd }\\
\quad \ldots & \quad \ldots \quad\ldots \\
s_{i_{2k+2}} = s_{i_{2k}}, &\quad \mbox{if}\; m_{s_{2k-1}s_{2k}}\; \mbox{is even}\\
s_{i_{2k+2}} = s_{i_{2k-1}}, &\quad \mbox{if}\; m_{s_{2k-1}s_{2k}}
\; \mbox{is odd}.
\end{cases}
$$
In \cite{Bokut-Shiao}, It was pointed out that the hypothesis holds
for Groebner-Shirshov bases for finite Coxeter groups. In the
following we point out that the hypothesis not necessarily correct.
Now consider the Coxeter matrix $A_4=\begin{bmatrix}
1 & 3 & 2 & 2 \\
3 & 1 & 3 & 2 \\
2 & 3 & 1 & 3 \\
2 & 2 & 3 & 1
\end{bmatrix}$.
By using above notation $s_is_i=1$ for $i=1,2,3,4$ we have
\begin{align}
&m(s_2,s_1)=s_2s_1s_2=m(s_1,s_2)=s_1s_2s_1\quad\mbox{since}\quad m_{21}=m_{12}=3\notag\\
&m(s_3,s_2)=s_3s_2s_3=m(s_2,s_3)=s_2s_3s_2\quad\mbox{since}\quad m_{32}=m_{23}=3\notag\\
&m(s_3,s_1)=s_3s_1=m(s_1,s_3)=s_1s_3\quad \mbox{since}\quad
m_{31}=m_{13}=2\notag \ {\rm and} \\
&m(s_4,s_3)=s_4s_3s_4=m(s_3,s_4)=s_3s_4s_3\notag\\
&m(s_4,s_2)=s_4s_2=m(s_2,s_4)=s_2s_4\notag\\
&m(s_4,s_1,)=s_4s_1=m(s_1,s_4)=s_1s_4\notag.
\end{align}
These are all initial relations. Now one can notice that there is
extra elements for $A_4$ such as:
$$(m-1)(s_3,s_2)m(s_1,s_3)=m(s_2,s_3)(m-1)(s_1,s_3)$$ which gives
the relation $s_3s_2s_1s_3=s_2s_3s_2s_1$, and similarly,
$(m-1)(s_4,s_3)(m-1)(s_2,s_4)m(s_1,s_4)=m(s_3,s_4)(m-1)(s_2,s_4)(m-1)(s_1,s_4)$
which gives the $s_4s_3s_2s_1s_4=s_3s_4s_3s_2s_1$. Thus the second
element $(s_2,s_4)$ shows the complexity of this notation in the
article \cite{Bokut-Shiao}. Thus $(m-1)(s_4,s_3) m(s_2,s_4) =
m(s_3,s_4)(m-1)(s_2,s_4)$  that gives $s_4s_3s_2s_4=s_4s_4s_3s_2$.

Now note that if we take the previous pair as an index the second
component of a pair can be determined by being even or odd of the
associated entry of the Coxeter matrix.  But the first component of
the chosen can be any element less than the element in the second
component. Thus if we write the pairs procedure in the article, the
second element of the pair $(s_1,s_2)$ must be $s_2$ but it says
that the element must be $s$ or $s'$.

%YANÝ BÝR ÝKÝLÝNÝN ÝKÝNCÝ ELEMANI KENDÝNDEN ONCEKÝ ÝKÝLYÝ ÝNDEX
%OLARAK ALDIÐIMIZDA MATRÝXTE KARÞIGELEN DEGERÝN TEK VEYA CÝFT OLMASI
%ÝLE BELÝRLENÝYOR. iKLÝNÝN ÝLK ELEMANI ÝSE ÝKÝNCÝ ELEMANDAN KÜÇÜK
%HERHANGÝ BÝR ELEMAN OLABÝLÝR. MAKALEDEKÝ GÝBÝ YAZARSAK SANKÝ HEP
%ÝKÝNCÝ ELEMAN $(s_1,s_2)$ OLAMLI GÝBÝ AMA KENDÝSÝDE $s_2$ NÝN $s_2$
%DEGÝL $s$ VEYA $s'$ OLDUGUNU SOYLUYOR.

\subsection{A counter-example to Hypothesis.}

The Coxeter matrix of affine Weyl group ${\widetilde{A}}_3$ is
$M=\begin{bmatrix}
1 & 3 & 2 & 2 \\
3 & 1 & 3 & 2 \\
2 & 3 & 1 & 3 \\
2 & 2 & 3 & 1
\end{bmatrix}$.

\noindent Initial relations are

$(s_0)^2=(s_1)^2=(s_2)^2=(s_3)^2=1$,

$s_0s_1s_0=s_1s_0s_1$,

$s_0s_2=s_2s_0$,

$s_0s_3s_0=s_3s_0s_3$,

$s_1s_2=s_1s_2s_1s_2$,

$s_1s_3=s_3s_1$,

$s_2s_3s_2=s_3s_2s_3$.

$$
\,
$$

\noindent The Groebner-Shirshov Basis for $W$ consists of the
elements

$s_0s_0-1$,

$s_1s_1-1$,

$s_2s_2-1$,

$4s_3s_3-1$,

$s_1s_0s_1-s_0s_1s_0$,

$s_2 s_1 s_2 - s_1 s_2 s_1$,

$s_3 s_2 s_3 - s_2 s_3 s_2$,

$s_3 s_0 s_3 - s_0 s_3 s_0$,

$s_3 s_1 - s_1 s_3$,

$s_2 s_0 - s_0 s_2$

$s_3 s_0 s_1 s_0 - s_1 s_3 s_0 s_1$,

$s_3 s_0 s_2 s_3 s_2 - s_0 s_3 s_0 s_2 s_3$,

$s_3 s_0 s_2 s_3 s_0 - s_2 s_3 * s_0 s_2 s_3$,

$s_3 s_2 s_1 s_3 - s_2 s_3 s_2 s_1$,

$s_3 s_0 s_1 s_3 - s_0 s_3 s_0 s_1$,

$s_2 s_1 s_0 s_2 - s_1 s_2 s_1 s_0$,

$s_3 s_0 s_2 s_1 s_3 s_2 s_1 - s_0 s_3 s_0 s_2 s_1 s_3 s_2$,

$s_3 s_0 s_1 s_2 s_3 s_2 - s_0 s_3 s_0 s_1 s_2 s_3$,

$s_3 s_2 s_1 s_0 s_3 s_0 - s_2 s_3 s_2 s_1 s_0 s_3$,

$s_3 s_0 s_2 s_1 s_3 s_0 s_1 - s_2 s_3 s_0 s_2 s_1 s_3 s_0$,

$s_3 s_0 s_1 s_2 s_3 s_0 s_2 - s_0 s_3 s_0 s_1 s_2 s_3 s_0$,

$s_3 s_0 s_2 s_1 s_3 s_0 s_2 s_1 s_0 - s_0 s_3 s_0 s_2 s_1 s_3 s_0
s_2 s_1$,

$s_ 3 s_0 s_1 s_2 s_1 s_3 s_2 s_1 - s_0 s_3 s_0 s_1 s_2 s_1 s_3
s_2$,

$s_3 s_0 s_2 s_1 s_0 s_3 s_0 s_1 - s_2 s_3 s_0 s_2 s_1 s_0 s_3 s_0$,

$s_3 s_0 s_1 s_2 s_1 s_3 s_0 s_2 s_1 s_0 - s_0 s_3 s_0 s_1 s_2 s_1
s_3 s_0 s_2 s_1$,

$s_3 s_0 s_1 s_2 s_3 s_0 s_1 s_2 s_1 - s_0 s_3 s_0 s_1 s_2 s_3 s_0
s_1 s_2$,

$s_3 s_0 s_1 s_2 s_1 s_3 s_0 s_1 s_2 s_1 s_0 - s_0 s_3 s_0 s_1 s_2
s_1 s_3 s_0 s_1 s_2 s_1$.

Next let us check if these words are appropriate for the hypothesis
or not.

$s_3s_0s_1s_0-s_1s_3s_0s_1=(m-1)(s_3,s_1)m(s_0,s_1)-m(s_1,s_3)$,

$s_3s_0s_2s_3s_2-s_0s_3s_0s_2s_3=(m-1)(s_3,s_0)m(s_2,s_3)-m(s_0,s_3)(m-1)(s_2,s_3)$,

$s_3 s_0 s_2 s_3 s_0 - s_2 s_3 * s_0 s_2
s_3=(m-1)(s_3,s_2)m(s_0,s_3)-m(s_2,s_3)(m-1)(s_0,s_3)$,

$s_3 s_2 s_1 s_3 - s_2 s_3 s_2
s_1=(m-1)(s_3,s_2)m(s_0,s_3)-m(s_2,s_3)(m-1)(s_0,s_3)$,

$s_3 s_0 s_1 s_3 - s_0 s_3 s_0
s_1=(m-1)(s_3,s_0)m(s_1,s_3)-m(s_0,s_3)(m-1)(s_1,s_3)$,

$s_2 s_1 s_0 s_2 - s_1 s_2 s_1
s_0=(m-1)(s_2,s_1)m(s_0,s_2)-m(s_1,s_2)(m-1)(s_0,s_2)$,

$ s_3s_0s_2s_1s_3s_2s_1-s_0s_3s_0s_2s_1s_3s_2\equiv
(m-1)(s_3,s_0)(m-1)(s_2,s_3)m(s_1,s_2)
-m(s_0,s_3)(m-1)m(s_2,s_3)(m-1)(s_1,s_2)$,

$s_3s_0s_1s_2s_3s_2-s_0s_3s_0s_1s_2s_3\equiv(m-1)(s_3,s_2)(m-1)(s_1,s_3)m(s_2,s_3)
-m(s_2,s_3)(m-1)(s_1,s_3)(m-1)(s_2,s_3)$,

$s_3s_2s_1s_0s_3s_0-s_2s_3s_2s_1s_0s_3\equiv(m-1)(s_3,s_2)(m-1)(s_1,s_3)m(s_0,s_3)
-m(s_2,s_3)(m-1)(s_1,s_3)(m-1)(s_0,s_3)$.

Thus we see that the elements above in the Groebner-Shirshov basis
can be written as the procedure in the hypothesis. But we have a
problem for the following element.

Let consider $s_3 s_0 s_2 s_1 s_3 s_0 s_1 - s_2 s_3 s_0 s_2 s_1 s_3
s_0$. Since the second term of the element begins with $s_2 s_3$, we
must start by $(m-1)(s_3,s_2)=s_3 s_2$.

Since $m_{43}=3$($s_3$ and $s_2$ are 4th and 3rd elements
respectively) is an odd number, we should continue by $m(a,s_0)$
where $a$ is a reflection less than  $s_0$ however this is not
allowed in the hypothesis. Actually this term can be written as
$(m-1)(s_3,s_2)(m-1)(s_0,s_3)m(s_1,s_0)-m(s_2,s_3)(m-1)(s_0,s_3)$

$(m-1)(s_1,s_0),$ but it does not follow the procedure in the
hypothesis. Similarly the following elements can be written as the
procedure in the hypothesis.

$s_3 s_0 s_1 s_2 s_3 s_0 s_2 - s_0 s_3 s_0 s_2 s_1 s_3 s_0 s_2 s_1
\equiv (m-1)(s_3,s_0)(m-1)(s_1,s_3)(m-1)(s_2,s_3)m(s_0,s_2)
-m(s_0,s_3)(m-1)(s_1,s_3)(m-1)(s_2,s_3)(m-1)(s_0,s_2)$,

$s_3s_0s_2s_1s_3s_0s_2s_1s_0-s_0s_3s_0s_2s_1s_3s_0s_2s_1s_1 \equiv
(m-1)(s_3,s_0)(m-1)(s_2,s_3)(m-1)(s_1,s_2)m(s_0,s_1)
-m(s_0,s_3)(m-1)(s_2,s_3)(m-1)(s_1,s_2)(m-1)(s_0,s_1)$,

$s_3 s_0 s_1 s_2 s_1 s_3 s_2 s_1 - s_0 s_3 s_0 s_1 s_2 s_1 s_3 s_2
\equiv (m-1)(s_3,s_0)(m-1)(s_1,s_3)(m-1)(s_2,s_3)m(s_1,s_2)
-m(s_0,s_3)(m-1)(s_1,s_3)(m-1)(s_2,s_3)(m-1)(s_1,s_2)$.

\noindent The element
$s_3s_0s_2s_1s_0s_3s_0s_1-s_2s_3s_0s_2s_1s_0s_3s_0$ is written as
$(m-1)(s_3,s_2)(m-1)(s_1,s_3)(m-1)(s_0,s_3)m(s_1,s_0)-m(s_3,s_2)(m-1)(s_1,s_3)(m-1)(s_0,s_3)(m-1)(s_1,s_0),$
but the last term of the expression does not obey the procedure of
the hypothesis.\\

\noindent The remaining elements obey the rule in the hypothesis.

$s_3s_0s_1s_2s_1s_3s_0s_2s_1s_0 - s_0s_3s_0s_1s_2s_1s_3s_0r-2s_1
\equiv
(m-1)(s_3,s_0)(m-1)(s_1,s_3)(m-1)(s_2,s_3)(m-1)(s_1,s_2)m(s_0,s_1)
-m(s_0,s_3)(m-1)(s_1,s_3)(m-1)(s_2,s_3)(m-1)(s_1,s_2)(m-1)(s_0,s_1)$,

$s_3s_0s_1s_2s_3s_0s_1s_2s_1 - s_0s_3s_0s_1s_2s_3s_0s_1s_2 \equiv
(m-1)(s_3,s_0)(m-1)(s_1,s_3)(m-1)(s_2,s_3)(m-1)(s_0,s_2)m(s_1,s_2)
-m(s_0,s_3)(m-1)(s_1,s_3)(m-1)(s_2,s_3)(m-1)(s_0,s_2)(m-1)(s_1,s_2)$,

$s_3s_0s_1s_2s_1s_3s_0s_1s_2s_1s_0-s_0s_3s_0s_1s_2s_1s_3s_0s_1s_2s_1
\equiv
(m-1)(s_3,s_0)(m-1)(s_1,s_3)(m-1)(s_2,s_3)(m-1)(s_0,s_2)(m-1)(s_1,s_2)m(s_0,s_1)
-m(s_0,s_3)(m-1)(s_1,s_3)(m-1)(s_2,s_3)(m-1)(s_0,s_2)(m-1)(s_1,s_2)m(s_0,s_1)$.

\noindent So we conclude that two elements can not be expressed as
the procedure in the hypothesis. With this counter-example, we
obtain that the hypothesis does not hold for positive degenerate
affine Weyl groups.

\section{Calculation of Groebner-Shirshov bases for infinite Weyl Groups.}

In this section our aim is to calculate Grobner-Shirshov basis of
positive degenerate infinite affine Weyl group $\widetilde{A}_n$
which is isomorphic to semi-direct product of permutation group
$S_n$ with translation group $\mathbb Z^{n-1}$, where $\mathbb
Z^{n-1}$ is the eigen-module of the action of the permutation group
on its representation module $\mathbb Z^{n}$. It has a presentation
with generators $\{s_0,s_1,\ldots,s_n\}$ and relations

\begin{align}
s_is_i&=1,\, i=0,\ldots,n;\notag\\
s_i s_j&=s_j s_i, \,i=0,\ldots,n-2,j-i>1;\notag\\
s_i s_{i+1}s_i&=s_{i+1}s_is_{i+1},\, i=0,\ldots,n-1 \,\text{and}\notag\\
s_0s_ns_0&=s_ns_0s_n.\notag
\end{align}

\noindent Here we will choose the linear ordering
$s_0>s_1>s_2>\cdots
>s_n$ on the generators on the contrary of the ordering chosen by
Bokut and Shiao in \cite{Bokut-Shiao}.

\begin{thm}\label{civan-arda}
Let $\widetilde{A_n}$ be affine Weyl group generated by
$s_0,s_1,s_2,\ldots , s_n$ with defining relations

$s_i^2=1, i=0,\ldots,n$

$s_i s_j=s_j s_i, i=0,\ldots,n-2, j=2,\ldots,n, j-i>1$

$s_i s_{i+1} s_i= s_{i+1}s_i s_{i+1}, i=0,\ldots,n-1$ and

$s_0 s_n s_0=s_n s_0 s_n$.

\noindent Let us define the words
$$s_{ij}=\left\{%
\begin{array}{lll}
    s_i s_{i-1} \cdots s_{j}, & i>j; \\
    s_i                    ,& i=j;\\
    s_i s_{i+1} \cdots s_j, & i<j. \\
\end{array}%
\right.    $$ and $\widehat{s}_{ij}=s_i s_{i-1} s_{i+1}s_i \cdots
s_{j+1} s_j$ where $j \geq i-1$.

If we identify a relation $u=v$ with polynomial $u-v$, a
Gr\"{o}bner-Shirshov basis for $\widetilde{A_n}$ with respect to
Deglex order with $s_0 > s_1 > \cdots > s_n$ consists of initial
relations together with the following polynomials.

$(1)$ $s_{ij}s_i - s_{i+1}s_{ij}$ where $j>i, i=0,\ldots,n-2,
j=i+2,\ldots,n$ when $i = 0$ and $j\neq n$.

$(2)$ $s_0 s_{nk} s_j - s_j s_0 s_{nk}$ where $j=2,\ldots,n-2$, $k =
n, n-1, \ldots, j+2$.

$(3)$ $s_0 s_{nj} s_{j+1}  - s_{j} s_0 s_{nj}$ where $j = n-1,
\ldots, 0$.

$(4)$ $s_0 s_{nj} s_{0}  - s_{n} s_0 s_{nj}$ where $j = 2, \ldots,
n-1$.

$(5)$ $s_0 s_{nj}\widehat{s }_{1k}s_{k+1}- s_n s_0
s_{nj}\widehat{s}_{1k}$ where $j=2,\ldots,n-1,n$, and
$k=0,\ldots,n-1$.

$(6)$  $s_{0j} s_{n} s_{0} s_n  - s_{1} s_{0j} s_{n} s_{0}$ where $j
= 1, \ldots, n-1$.

$(7)$ $s_0 s_{nj} s_1 s_0 s_{nk} s_1 - s_n s_0 s_{nj}s_1 s_0 s_{nk}$
where $j=2,\ldots,n-1$, and  $k = j-1,\ldots,n$.

$(8)$ $s_0 s_{nj} s_1 s_0 s_{nk} \widehat{s}_{2l}s_{l+1} - s_n s_0
s_{nj}s_1 s_0 s_{nk}\widehat{s}_{2l}$ where $j=2,\ldots,n-1$,
$k=j+1,\ldots,n$, \text{and}  $l = 1,\ldots, n-1$.

$(9)$ $s_{0j} s_n s_k s_0 s_{nk}- s_1 s_{0j}s_n s_0 s_{n-1,k}
s_{k+1}$ where $j=1,\ldots,n-1$, and $k = n-1, n-2, \ldots, j + 2$.

$(10)$  $s_{0j}s_n s_{kl}s_0 s_{nl} - s_1 s_{0j}s_n s_0
s_{n-1,l}s_{k+1,l+1}$ where $j=1,\ldots,n-1$, $k = n-1,
n-2,\ldots,j+2$ and $l = k-1,\ldots,2$.
\end{thm}
\begin{proof}
The proof constitute two parts. In Part I we investigate where the
new binomials come from.

$(1)$ $s_{ij}s_i - s_{i+1}s_{ij}$ where $j>i, i=0,\ldots,n-2,
j=i+2,\ldots,n$ when $i = 0$ and $j\neq n$.

For $j = i + 2$ we have
\begin{align}
\langle s_i s_{i+1}s_i - s_{i+1}s_{i} s_{i+1}; s_{i} s_{i+2} -
s_{i+2} s_{i}\rangle & = s_i s_{i+1}s_{i+2}s_i - s_{i+1}s_{i}
s_{i+1} s_{i+2}\notag\\
&= s_{i,i+2}s_i - s_{i+1}s_{i,i+2}\notag.
\end{align}

For $j=i+3,\ldots,n$ the others are obtained by induction on $j$.
$$
\langle s_{i,j-1}s_i - s_{i+1}s_{i,j-1}; s_{i}s_j -
s_{j}s_{i}\rangle = s_{ij}s_i - s_{i+1}s_{ij}.
$$

$(2)$ $s_0 s_{nk} s_j - s_j s_0 s_{nk}$ where $j=2,\ldots,n-2$, $k =
n, n-1, \ldots, j+2$.

For $k = n$ and $j=2,\ldots,n-2$ we have
$$
\langle s_0  s_j - s_j s_0 ; s_j s_{n} - s_n s_{j}\rangle = s_0
s_{n} s_j - s_j s_0 s_{n}.
$$
For $k = n-1,\ldots,j+2$ the others are obtained by induction on
$k$.
$$
\langle s_0 s_{n,k+1} s_j - s_j s_0 s_{n,k+1}; s_j s_{k} - s_k
s_{j}\rangle = s_0 s_{nk} s_j - s_j s_0 s_{nk}.
$$

$(3)$ $s_0 s_{nj} s_{j+1}  - s_{j} s_0 s_{nj}$ where $j = n-1,
\ldots, 0$.

For $j = n-1$ we have
$$
\langle s_0  s_{n-1} - s_{n-1} s_0 ;  s_{n-1} s_{n}s_{n-1}- s_n
s_{n-1}s_{n}\rangle = s_0 s_{n} s_{n-1}s_n - s_{n-1} s_0
s_{n}s_{n-1}.
$$
The others are obtained by composition of the element of the type
$(2)$ taking $k=j+2$.
$$
\langle s_0 s_{n,j+2} s_j - s_j s_0 s_{n,j+2}; s_{j} s_{j+1}s_{j}-
s_{j+1} s_{j}s_{j+1}\rangle = s_0 s_{nj} s_{j+1}  - s_{j} s_0
s_{nj}.
$$

$(4)$ $s_0 s_{nj} s_{0}  - s_{n} s_0 s_{nj}$ where $j = 2, \ldots,
n-1$.

For $j = n-1$ we have
$$
\langle s_0  s_{n} s_0 - s_{n} s_0 s_{n} ; s_0 s_{n-1} -
s_{n-1}s_{0}\rangle = s_0 s_{n} s_{n-1}s_0 - s_{n} s_0 s_{n}s_{n-1}.
$$

For $j = n-2,\ldots,2$ the others are obtained by induction on $j$.
$$
\langle s_0 s_{n,j+1} s_0 - s_n s_0 s_{n,j+1}; s_{0} s_{j} - s_{j}
s_{0}\rangle = s_0 s_{nj} s_{0}  - s_{n} s_0 s_{nj}.
$$

$(5)$ $s_0 s_{nj}\widehat{s }_{1k}s_{k+1}- s_n s_0
s_{nj}\widehat{s}_{1k}$ where $j=2,\ldots,n-1,n$, and
$k=0,\ldots,n-1$.

For $j = n$ and $k = 0$  we have
$$
\langle s_0  s_{n} s_0 - s_{n} s_0 s_{n} ; s_0 s_{1} s_0  - s_{1}
s_{0} s_{1} \rangle = s_0 s_{n} s_{1}s_0 s_1 - s_{n} s_0
s_{n}s_{1}s_0.
$$

For $j = n-1,\ldots,2$ and $k = 0$  we have
$$
\langle s_0  s_{nj} s_0 - s_{n} s_0 s_{nj} ; s_0 s_{1} s_0  - s_{1}
s_{0} s_{1} \rangle = s_0 s_{nj} s_{1} s_0 s_1 - s_{n} s_0
s_{nj}s_{1}s_0.
$$
The others are obtained by induction on $k$.

$$
\langle s_0 s_{nj}\widehat{s }_{1,k-1}s_{k}- s_n s_0
s_{nj}\widehat{s}_{1,k-1}; s_{k} s_{k+1}s_{k}- s_{k+1}
s_{k}s_{k+1}\rangle = s_0 s_{nj}\widehat{s }_{1k}s_{k+1}- s_n s_0
s_{nj}\widehat{s}_{1k}.
$$

$(6)$  $s_{0j} s_{n} s_{0} s_n  - s_{1} s_{0j} s_{n} s_{0}$ where $j
= 1, \ldots, n-1$.

For $j=1$ we have

$$
\langle  s_0 s_{1} s_0  - s_{1} s_{0} s_{1}; s_0  s_{n} s_0 - s_{n}
s_0 s_{n} \rangle = s_0 s_{1} s_{n} s_0 s_n - s_{1} s_0 s_{1}
s_{n}s_0.
$$

The others are obtained by composition of the element of the type
$(1)$ taking $i=0$.

$$
\langle  s_{0j} s_0  - s_{1} s_{0j} ; s_0  s_{n} s_0 - s_{n} s_0
s_{n} \rangle = s_{0j} s_{n} s_{0} s_n  - s_{1} s_{0j} s_{n} s_{0}
$$
where $j = 2,\ldots, n-1$.

$(7)$ $s_0 s_{nj} s_1 s_0 s_{nk} s_1 - s_n s_0 s_{nj}s_1 s_0 s_{nk}$
where $j=2,\ldots,n-1$, and  $k = j-1,\ldots,n$.

For $k=n$ it is obtained by composition of the element of the type
$(5)$ taking $k=0$.

$$
\langle  s_0 s_{nj} s_1 s_0 s_1 - s_{n}s_0 s_{nj} s_1 s_0 ; s_1
s_{n} - s_{n} s_1  \rangle = s_0 s_{nj} s_1 s_0 s_{n} s_1 - s_n s_0
s_{nj}s_1 s_0 s_{n}.
$$
Note that we do not consider the case where $j = n$ because one can
show that it reduces to zero.

The others are obtained by induction on $k$.
$$
\langle s_0 s_{nj} s_1 s_0 s_{n,k+1} s_1 - s_n s_0 s_{nj}s_1 s_0
s_{n,k+1}; s_1 s_{k} - s_{k} s_1  \rangle = s_0 s_{nj} s_1 s_0
s_{nk} s_1 - s_n s_0 s_{nj}s_1 s_0 s_{nk}
$$
where $k = n-1,\ldots, j-1$.

$(8)$ $s_0 s_{nj} s_1 s_0 s_{nk} \widehat{s}_{2l}s_{l+1} - s_n s_0
s_{nj}s_1 s_0 s_{nk}\widehat{s}_{2l}$ where $j=2,\ldots,n-1$,
$k=j+1,\ldots,n$, \text{and}  $l = 1,\ldots, n-1$.

For $l = 1$ it is obtained by composition of the element of the type
$(7)$.
$$
\langle s_0 s_{nj} s_1 s_0 s_{nk} s_1 - s_n s_0 s_{nj}s_1 s_0
s_{nk}; s_1 s_{2} s_1 - s_{2} s_1 s_{2} \rangle = s_0 s_{nj} s_1 s_0
s_{nk} s_{2} s_1 s_{2} - s_n s_0 s_{nj}s_1 s_0 s_{nk}s_{2}s_{1}.
$$
The furthers are obtained by using the induction on $l$.
$$
\langle s_0 s_{nj} s_1 s_0 s_{nk} \widehat{s}_{2,l-1}s_{l} - s_n s_0
s_{nj}s_1 s_0 s_{nk}\widehat{s}_{2,l-1}; s_l s_{l+1} s_l - s_{l+1}
s_l s_{l+1} \rangle = s_0 s_{nj} s_1 s_0 s_{nk}
\widehat{s}_{2l}s_{l+1} - s_n s_0 s_{nj}s_1 s_0
s_{nk}\widehat{s}_{2l}.
$$

$(9)$ $s_{0j} s_n s_k s_0 s_{nk}- s_1 s_{0j}s_n s_0 s_{n-1,k}
s_{k+1}$ where $j=1,\ldots,n-1$, and $k = n-1, n-2, \ldots, j + 2$.

\noindent It is obtained by composition of the elements of the types
$(6)$ and $(3)$ respectively.
$$
\langle s_{0j} s_{n} s_{0} s_n  - s_{1} s_{0j} s_{n} s_{0}; s_0
s_{nk} s_{k+1}  - s_{k} s_0 s_{nk} \rangle = s_{0j} s_n s_k s_0
s_{nk}- s_1 s_{0j}s_n s_0 s_{n-1,k} s_{k+1}.
$$

$(10)$  $s_{0j}s_n s_{kl}s_0 s_{nl} - s_1 s_{0j}s_n s_0
s_{n-1,l}s_{k+1,l+1}$ where $j=1,\ldots,n-1$, $k = n-1,
n-2,\ldots,j+2$ and $l = k-1,\ldots,2$.

\noindent For $l = k-1$, it is obtained by composition of the
elements of the types $(9)$ and $(3)$ respectively.

\begin{align}
\langle s_{0j} s_n s_k s_0 s_{nk} - s_1 s_{0j}s_n s_0 s_{n-1,k}
s_{k+1} ; s_0 s_{n,k-1} s_{k}  - s_{k-1} s_0 s_{n,k-1}\rangle & =
s_{0j} s_n s_k s_{k-1}s_0 s_{n,k-1} - s_1 s_{0j}s_n
s_{n-1,k}\notag\\
& s_{k+1} s_{k-1} s_{k}\notag\\
&= s_{0j} s_n s_k s_{k-1}s_0 s_{n,k-1} - s_1 s_{0j}s_n
s_{n-1,k-1}\notag\\& s_{k+1} s_{k}\notag.
\end{align}
The others are obtained by induction on $l$.
\begin{align}
\langle s_{0j}s_n s_{k,l+1}s_0 s_{n,l+1} - s_1 s_{0j}s_n s_0
s_{n-1,l+1}s_{k+1,l+2} ; s_0 s_{nl} s_{l+1} - s_{l} s_0
s_{nl}\rangle & = s_{0j}s_n s_{k,l+1} s_l s_0 s_{nl} - s_1 s_{0j}s_n
s_0 s_{n-1,l+1} \notag\\
& s_{k+1,l+2} s_{l} s_{l+1} \notag\\
&= s_{0j}s_n s_{kl}s_0 s_{nl} - s_1 s_{0j}s_n s_0
s_{n-1,l}s_{k+1,l+1}\notag.
\end{align}

Part II : Reduction.

 Now we have to show that all other compositions
are trivial (reduces to zero) relative to initial relations
including Groebner-Shirshov basis. Notice that a composition
$\langle f;g \rangle$ is trivial if leading words of $f$ and $g$ has
no overlap. Because of that, we will only look compositions with
overlap leading words.

In this part of the proof we have a lemma which will help  us to
reduce some compositions to zero easily.

\begin{lem}\label{abbas}
Let $f = f_1 - f_2$, $g = g_1 - g_2$ and $h = h_1 - h_2$ be
binomials. Assume that  $\overlap (f_1 , g_1) = v_1$, that is $f_1 =
a_1 v_1$ and $g_1 = v_1 b_1$, and $\overlap (g_1 , h_1) = v_2$, that
is $g_1 = a_2 v_2$ and $h_1 = v_2 b_2$. Define
$$
\langle f;g \rangle = (f_1 - f_2) b_1 - a_1 (g_1 - g_2) = a_1 g_2 -
f_2 b_1
$$
and

$$
\langle g;h \rangle = (g_1 - g_2) b_2 - a_2 (h_1 - h_2) = a_2 h_2 -
g_2 b_2.
$$
If $\langle f;g \rangle$ is trivial and $a_2 = v_1 \bar{a}_2$, then
$\langle f; a_2 h_2 - g_2 b_2 \rangle$ is also trivial.
\end{lem}
\begin{proof}
\begin{align}
\langle f; a_2 h_2 - g_2 b_2 \rangle & = (f_1 - f_2)\bar{a}_2 h_2 -
a_1 (a_2 h_2 - g_2 b_2)\notag\\
&= a_1 g_2 b_2 - f_2 \bar{a}_2 h_2\notag\\
& = a_1 g_2 b_2 + f_2 \bar{a}_2(h_1 - h_2) -  f_2 \bar{a}_2 h_1
\notag.
\end{align}
Since $h_1 = v_2 b_2$ and $b_1 = \bar{a}_2 h_1$, we have
\begin{align}
\langle f; a_2 h_2 - g_2 b_2 \rangle & = a_1 g_2 b_2 + f_2 \bar{a}_2(h_1 - h_2) -  f_2 b_1 b_2\notag\\
&= (a_1 g_2 -f_2 b_1) b_2 + f_2 \bar{a}_2(h_1 - h_2)\notag\\
& = 0\notag.
\end{align}
\end{proof}
Now we can start to show that the remaining composition reduces to
zero. First of all we look the composition where $s_0 s_1 s_0 - s_1
s_0 s_1$ is the left element. Since almost all element types start
with $s_0$, we have a lot compositions to check. However, by Lemma
\ref{abbas} we will only show three such composition reducing to
zero. The others will also reduce to zero as an application of Lemma
\ref{abbas}. So let us check  three such compositions.
\begin{align}
\langle s_0 s_1 s_0 - s_1 s_0 s_1 ; s_0 s_1 s_2 s_0 - s_1 s_0 s_1
s_2\rangle & = s_0 s_1 s_1 s_0 s_1 s_2 - s_1 s_0 s_1 s_1 s_2
s_0\notag\\
& = s_0 (s_1 s_1 -1) s_0 s_1 s_2 - s_1 s_0 (s_1 s_1 - 1) s_2 s_0
\notag\\ &+ (s_0 s_0 - 1) s_1 s_2 - s_1 (s_0 s_2 - s_2 s_0) s_0 -
s_1 s_2 (s_0
s_0 - 1)\notag\\
& = 0.\notag
\end{align}
$$
\langle s_0 s_1 s_0 - s_1 s_0 s_1 ; s_0 s_j  - s_j s_0 \rangle = s_0
s_1 s_j s_0 - s_1 s_0 s_1 s_j
$$
Notice that for $j=2$ this is relation with type $(1)$ choosing $i =
0, j=2$. If $j\neq 2$ then
\begin{align}
s_0 s_1 s_j s_0 - s_1 s_0 s_1 s_j & = s_0 (s_1 s_j - s_j s_1) s_0 -
s_1 s_0 (s_1 s_j - s_j s_1) + (s_0 s_j - s_j s_0)s_1 s_0 \notag\\ &-
s_1 (s_0 s_j - s_j s_0) s_1  -  (s_1 s_j - s_j s_1) s_0 s_1 + s_j
(s_0 s_1 s_0 - s_1 s_0 s_1)\notag.
\end{align}
We also have

$$
\langle s_0 s_1 s_0 - s_1 s_0 s_1 ; s_0 s_n s_0  - s_n s_0
s_n\rangle = s_0 s_n s_1 s_0 s_1 - s_n s_0 s_n s_1 s_0
$$
which is relation with type $(6)$ choosing $j=1$.

The other compositions with left element $s_0 s_1 s_0 - s_1 s_0 s_1$
reduce to zero as an application of Lemma \ref{abbas}.

Now consider compositions in which $s_0 s_n s_0  - s_n s_0 s_n$ is
left element. By above argument we only have to look the following
compositions.
\begin{align}
\langle s_0 s_n s_0 - s_n s_0 s_n ; s_0 s_1 s_2 s_0 - s_1 s_0 s_1
s_2\rangle & = s_0 s_n s_1 s_0 s_1 s_2 - s_n s_0 s_n s_1 s_2
s_0\notag\\
& = (s_0 s_n s_1 s_0 s_1 -  s_n s_0 s_n s_1 s_0 )s_2 + s_n s_0 s_n
s_1 (s_0 s_2 - s_2 s_0)\notag,
\end{align}
where $s_0 s_n s_1 s_0 s_1 -  s_n s_0 s_n s_1 s_0$ is relation with
type $(5)$ choosing $k = 0, j=n$.

We also have

$$
\langle s_0 s_n s_0 - s_n s_0 s_n ; s_0 s_j  - s_j s_0 \rangle = s_0
s_n s_j s_0  - s_n s_0 s_n s_j.
$$
Notice that for $j = n-1$ this is relation with type $(4)$ choosing
$j = n-2$. If $j\neq n-1$ then

\begin{align}
s_0 s_n s_j s_0 - s_n s_0 s_n s_j & = (s_0 s_n s_j - s_j s_0 s_n)
s_0
- s_n s_0 (s_n s_j - s_j s_n) - s_n (s_0 s_j - s_j s_0)s_n \notag\\
&+  (s_j s_n - s_n s_j) s_0 s_n  + s_j (s_0 s_n s_0 - s_n s_0
s_n)\notag.
\end{align}

The other compositions with left element $s_0 s_n s_0 - s_n s_0 s_n$
reduce to zero as an application of Lemma \ref{abbas}.

Now consider compositions with left element $s_{0j} s_0  - s_1
s_{0j}$ which is relation with type $(1)$ choosing $i = 0$.
\begin{align}
\langle s_{0j} s_0  - s_1 s_{0j} ; s_0 s_1 s_0 - s_1 s_0 s_1 \rangle
& = s_{0j} s_1 s_0 s_1 - s_1 s_{0j} s_1 s_0\notag\\
& = s_0 (s_{1j}  s_1  - s_2 s_{1j} ) s_0 s_1 - s_1 s_0 (s_{1j}  s_1  - s_2 s_{1j} ) s_0 \notag\\
& +  (s_0 s_2 - s_2 s_0) s_{1j} s_0 s_1  - s_1(s_0 s_2 - s_2 s_0) s_{1j} s_0\notag\\
& + s_2 (s_{0j}  s_1  - s_1 s_{0j} ) s_1 - s_1 s_2 (s_{0j}  s_0 -
s_1 s_{0j} ) \notag\\
& + s_2  s_1 s_0 (s_{1j}  s_1  - s_2 s_{1j} )  - (s_1 s_2 s_1 - s_2
s_1 s_2) s_{0j} \notag\\
& + s_2  s_1 (s_0 s_2 - s_2 s_0) s_{1j}\notag.
\end{align}
We have
$$
\langle s_{0j} s_0  - s_1 s_{0j} ; s_0 s_n s_0 - s_n s_0 s_n \rangle
= s_{0j} s_n s_0 s_n - s_1 s_{0j} s_n s_0
$$
which is relation with type $(6)$.

$$
\langle s_{0j} s_0  - s_1 s_{0j} ; s_0 s_k  - s_k s_0  \rangle =
s_{0j} s_k s_0  - s_1 s_{0j} s_k.
$$

We have the following cases to be analyzed:

Case $1)$ If $k = j$ then
$$
s_{0j} s_j s_0  - s_1 s_{0j} s_j  = s_{0,j-1} (s_j s_j -1) s_0 - s_1
s_{0,j-1} (s_j s_j -1) - (s_{0,j-1} s_0 - s_1 s_{0,j-1} = 0.
$$

Case $2)$ If $k = j + 1$ then it just gives $s_{0,j+1} s_0  - s_1
s_{0,j+1}$ which is relation with type $(1)$.

Case $3)$ If $k > j + 1$ then
$$
s_{0j} s_k s_0  - s_1 s_{0j} s_k = s_k ( s_{0j} s_0  - s_1 s_{0j})
$$
using commutativity of $s_k$.

Case $4)$ If $k < j $ then
\begin{align}
s_{0j} s_k s_0  - s_1 s_{0j} s_k
& = s_{0,k+1} s_k s_{k+2,j} s_0  -
s_1 s_{0,k+1} s_k s_{k+2,j}\quad(\text{by commutativity of
$s_k$})\notag\\
& = s_{0,k-1} (s_k s_{k+1}s_k - s_{k+1} s_{k} s_{k+1}) s_{k+2,j} s_0
- s_1 s_{0,k-1} (s_k s_{k+1}s_k - s_{k+1} s_{k} s_{k+1}) s_{k+2,j}
\notag\\
&+ s_{k+1} (s_{0j} s_0  - s_1 s_{0j})\notag.
\end{align}

The other compositions with left element $s_{0j} s_0  - s_1 s_{0j}$
reduce to zero as an application of Lemma \ref{abbas}.

Now consider compositions in which $s_0 s_{nj} s_0  - s_n s_0
s_{nj}$ is left element.

$$
\langle s_0 s_{nj} s_0  - s_n s_0 s_{nj} ; s_0 s_1 s_0 - s_1 s_0 s_1
\rangle = s_0 s_{nj} s_1 s_0 s_1 - s_n s_0 s_{nj} s_1 s_0
$$
which is relation with type $(5)$ choosing $k = 0, j =
n-1,\ldots,2$.

\begin{align}
\langle s_0 s_{nj} s_0  - s_n s_0 s_{nj} ; s_0 s_n s_0 - s_n s_0 s_n
\rangle &= s_0 s_{nj} s_n s_0 s_n - s_n s_0 s_{nj} s_n s_0\notag\\
& = s_0 s_n s_{n-1} s_{n} s_{n-2,j} s_0 s_n - s_n s_0 s_n s_{n-1}
s_{n-2,j} s_0 \quad(\text{by commutativity of $s_n$}) \notag\\
& = (s_0 s_n s_{n-1} s_{n}- s_{n-1}s_0 s_n s_{n-1}) s_{n-2,j} s_0
s_n  - s_n( s_0 s_n s_{n-1}s_n- \notag\\& s_{n-1}s_0 s_n s_{n-1})
s_{n-2,j}
s_0\notag\\
& + s_{n-1} (s_0 s_{nj} s_{0} - s_{n} s_0 s_{nj})  s_n  + s_{n-1}
s_n ( s_0 s_n s_{n-1}s_n- s_{n-1}s_0 s_n s_{n-1}) s_{n-2,j} \notag\\
& - s_n s_{n-1} (s_0 s_{nj} s_{0} - s_{n} s_0 s_{nj}) + (s_{n-1}
s_{n} s_{n-1} - s_{n} s_{n-1} s_{n})s_0 s_{nj}\notag,
\end{align}
where $s_0 s_n s_{n-1} s_{n}- s_{n-1}s_0 s_n s_{n-1}$ is relation
with type $(3)$ choosing $j = n-1$.

$$
\langle s_0 s_{nj} s_0  - s_n s_0 s_{nj} ; s_0 s_k  - s_k s_0
\rangle = s_0 s_{nj} s_k s_0  - s_n s_0 s_{nj} s_k.
$$
We have the following cases to be analyzed:

Case $1)$. Let $k < j-1$. Then
$$
s_0 s_{nj} s_k s_0  - s_n s_0 s_{nj} s_k  = (s_0 s_{n,j} s_k- s_k
s_0 s_{nj})s_0 + s_k (s_0 s_{n,j} s_0 - s_n s_0 s_{nj}) .
$$

Case $2)$. Let  $k = j - 1$. Then
$$
 s_0 s_{nj} s_{j-1} s_0  - s_n s_0 s_{nj} s_{j-1} = s_0 s_{n,j-1}s_0
 - s_n s_0 s_{n,j-1}.
$$

Case $3)$. Let $k = j$. Then
$$
s_0 s_{nj} s_j s_0  - s_n s_0 s_{nj} s_j = s_0 s_{n,j+1} ( s_{j} s_j
- 1)s_0 - s_n s_0 s_{n,j+1}( s_{j} s_j - 1) + (s_0 s_{n,j+1} s_0 -
s_n s_0 s_{n,j+1}).
$$

Case $4)$. Let $k > j $. Then
\begin{align}
s_0 s_{nj} s_k s_0  - s_n s_0 s_{nj} s_k & = s_0 s_{n,j+1}( s_j s_k
- s_k s_j)s_0 - s_n s_0 s_{n,j+1} ( s_j s_k - s_k s_j)\notag\\
& +   s_0 s_{n,j+1} s_k( s_j s_0 - s_0 s_j) + (s_0 s_{n,j+1} s_k s_0
- s_n s_0 s_{n,j+1}s_k) s_j\notag.
\end{align}

The other compositions with left element $s_{0} s_j  - s_j s_{0}$
reduce to zero as an application of Lemma \ref{abbas}.

\begin{align}
\langle s_{0} s_j  - s_j s_{0} ; s_j s_k  - s_k s_j  \rangle
& = s_{0} s_k s_j  - s_j s_{0} s_k \notag\\
& = (s_{0} s_k  - s_k s_{0})s_j - s_j(s_{0} s_k  - s_k s_{0}) \notag\\
& -  (s_j s_k  - s_k s_j)s_0 + s_k ( s_0 s_j  - s_j s_0).\notag
\end{align}
Similarly it can be shown that the following listed compositions
reduce to zero as an application of Lemma \ref{abbas}.

$$
\langle s_0 s_{n,k+1} s_j  - s_j s_0 s_{n,k+1} ; s_j s_t  - s_t s_j
s_1 \rangle \equiv 0 \Mod \, (\text{the listed Groebner-Shirshov
basis}) \quad \text{for}\quad t\neq k.
$$

$$
\langle s_0 s_{nk} s_j  - s_j s_0 s_{nk} ; s_j s_{j+1} s_j  -
s_{j+1} s_j s_{j+1} \rangle \equiv 0 \Mod \, (\text{the listed GS-
basis}) \quad \text{for}\quad  k\neq j+2.
$$

$$
\langle s_0 s_{nj} \hat{s}_{1k}s_{k+1}  - s_n s_0 s_{nj}\hat{s}_{1k}
; s_k s_{t}  - s_{t} s_k  \rangle \equiv 0 \Mod \, (\text{the listed
GS- basis}) \quad \text{for}\quad  k\neq 1, t \neq n.
$$

$$
\langle s_0 s_{nj} s_1 s_0 s_{nk} s_{1}  - s_n s_0 s_{nj}s_1 s_0 ;
s_1 s_{k}  - s_{k} s_1  \rangle \equiv 0 \Mod \, (\text{the listed
GS- basis}) \quad \text{for}\quad  k = j-2, j-3, \ldots,3.
$$

$$
\langle s_0 s_{nl} s_1 s_0 s_{nk}\hat{s}_{2l}s_{l+1}  - s_n s_0
s_{nj}s_1 s_0 s_{nk}\hat{s}_{2l} ; s_{l+1} s_{t}  - s_{t} s_{l+1}
\rangle \equiv 0 \Mod \, (\text{the listed GS- basis}).
$$

$$
\langle s_{0j} s_n s_0 s_{n}  - s_1 s_{0j} s_n s_0 ; s_0 s_{nk}
s_{k+1} - s_{k} s_0 s_{nk}  \rangle \equiv 0 \Mod \, (\text{the
listed GS- basis}) \quad \text{for}\quad  k = j+1, j, \ldots,1.
$$

$$
\langle s_{0j} s_n s_{k}s_0 s_{nk}  - s_1 s_{0j}s_0 s_{n-1,k}
s_{k+1} ; s_k s_{t}  - s_{t} s_k  \rangle \equiv 0 \Mod \,
(\text{the listed GS- basis}).
$$

$$
\langle s_{0j} s_n s_{k}s_0 s_{nk}  - s_1 s_{0j}s_0 s_{n-1,k}
s_{k+1} ; s_k s_{k+1} s_k - s_{k+1} s_k s_{k+1} \rangle \equiv 0
\Mod \, (\text{the listed GS- basis}).
$$

$$
\langle s_{0j} s_n s_{nk}s_0 s_{nl}  - s_1 s_{0j} s_n s_0 s_{n-1,l}
s_{k+1,l+1} ; s_l s_{l+1} s_l - s_{l+1} s_l s_{l+1} \rangle \equiv 0
\Mod \, (\text{the listed GS- basis}).
$$

$$
\langle s_{0j} s_n s_{nk}s_0 s_{nl}  - s_1 s_{0j} s_n s_0 s_{n-1,l}
s_{k+1,l+1} ; s_l s_{t}  - s_{t} s_l \rangle \equiv 0 \Mod \,
(\text{the listed GS- basis}).
$$
It completes the proof.

\end{proof}

\subsection{Representation of Reduced Elements in the group.}
The Composition Diamond Lemma tells that the irreducible elements
with respect to the Groebner-Shirsov basis give us a
representation for affine Weyl group $\widetilde{A}_n$. Since we
choose the total ordering $s_0 > s_1 > \cdots > s_n$ on the
generating set, the elements beginning with $s_0$ can not have
conjugate elements beginning with other generators. Therefore the
elements beginning with $s_0$ give us a representation of elements
with minimal length of in the flag subset $\widetilde{A}_n /
A_n$.\\

\noindent As a consequence of our calculations it is seen that
this representation can be expressed by words which are obtained
by juxtaposing of some blocks of the elements.

The first block in the row is a word satisfying the ordering
$\{s_0,s_1,\ldots s_n\}$ and the following conditions:

(i) $s_i s_j$ not in that block except $s_0s_n$, for $j-i>1,\,
i=0,\ldots,n-2$,

(ii) $s_0s_n s_j$ is not contained in the block for
$j=2,\ldots,n-2$, and

(iii) if the word is repeated in the sequent block again, the
conditions still hold.

By a simple calculation, we have these blocks in the following:

$s_0 s_n s_{n-1}s_1s_2\cdots s_{n-2}$,

$s_j s_{j+1} \cdots s_0s_ns_{n-1}s_1s_2 \cdots s_{j-1}$,

$s_1s_2 \cdots s_{n-2}s_0s_ns_{n-1}$,

$s_{n-1}s_1s_2 \cdots s_{n-2}s_0s_n$,

$s_ns_{n-1}s_1s_2 \cdots s_{n-2}s_0$,

$s_0s_ns_{n-1}s_{n-2} \cdots s_1$,

$s_js_{j-1} \cdots s_1s_0s_ns_{n-1}s_{n-2}\cdots s_{j+1}$,

$s_{n-1}s_{n-2} \cdots s_1 s_0 s_n$,

$s_n s_{n-1}s_{n-2} \cdots s_1 s_0$,

$s_{n-1}s_0s_{n}s_1s_2\cdots s_{n-2}$,

$s_j s_{j+1} \cdots s_{n-1}s_0s_{n}s_1s_2 \cdots s_{j-1}$,

$s_1s_2 \cdots s_{n-2}s_{n-1}s_0s_{n}$,

$s_{n}s_1s_2 \cdots s_{n-2}s_{n-1}s_0$,

$s_0s_{n}s_1s_2 \cdots s_{n-2}s_{n-1}$,

$s_{n-1}s_ns_{0}s_1s_2\cdots s_{n-2}$,

$s_j s_{j+1} \cdots s_{n-1}s_ns_{0}s_1s_2 \cdots s_{j-1}$,

$s_1s_2 \cdots s_{n-2}s_{n-1}s_ns_{0}$,

$s_{0}s_1s_2 \cdots s_{n-2}s_{n-1}s_n$,

$s_ns_{0}s_1s_2 \cdots s_{n-2}s_{n-1}$.

If we give the examples of this type words in the group
$\widetilde{A}_4$ the words beginning with $s_0$ are the
following:

$s_0s_1s_2s_3s_4,s_0s_4s_1s_2s_3,s_0s_4s_3s_1s_2,s_0s_4s_3s_2s_1$.

The other blocks satisfying the rules above can be easily written as
above.

Also there are blocks apart from this type blocks, with the
following properties:

(i*) they begin with $s_k$ and end with $s_k$ again,

(ii*) $3\leq \ell(\text{block}) \leq n$.

These blocks can be written as $s_k x s_k$ where $x$ is the string
of generators $\{s_0,\ldots,s_{k-1},s_{k+1},\ldots,s_n\}$ with
length $3\leq \ell(x) \leq n$.

If we give the example of this type words in the group
$\widetilde{A}_4$ the word beginning with $s_0$ is $s_0s_4s_1s_0$.

Of course there are a lot of blocks beginning with $s_1,s_2,s_3$ or
$s_4$. They can be written with respect to the rules above.

The representation of elements is made up of the juxtaposing of
obtained blocks. But when the blocks are juxtaposing, it is point
out that no leading terms of elements of Groebner-Shirshov basis
must be formed in the strings so that the obtained word can not be
reduced.

For example the block $(s_0 s_1 s_2 s_3 s_4)$ can take along no
block except itself.

%Ornegin $(s_0s_1s_2s_3s_4)$ blogu kendi haric yanina hicbir blok
%alamaz. Cunku $s_1s_2s_3s_4s_1-s_2s_1r-2s_3s_4$ Grobner shirshov
%tabaninin bir elemani olup $s_1s_2s_3s_4s_1$ elde etmemek icin
%blogumuz $s_1$ ile baslayan blok alamaz. Benzer nedenle
%$s_2s_3s_4s_2$ Grobner-Shirshov tabanindaki bir elemanin ilk
%(leading) kelimesi oldugundan $s_2$ ile baslayan blok alamaz. $s_3$
%le de baslayamaz cunku Grobner-Shirshov tabaninda
%$s_3s_4s_3-s_4s_3s_4$ var. Oyleyse $s_0$ baslayan blok alablir fakat
%$s_0s_1s_2s_3s_4s_0s_4$ yine Grobner-Shirshov tabanin bir elemaninin
%ilk (leading) kelime oldugundan $s_0$ dan sonra ancak $s_2$ gelir.
%Bu mantikla devam edilirse $(s_0s_1s_2s_3s_4)$ blogunu ancak kendisi
%takip edebildigi gorulur.

Therefore we have the following elements in the flag subset
$\widetilde{A}_4 / A_4$:

$(s_0s_1s_2s_3s_4)^p$,

$(s_0s_1s_2s_3s_4)^p s_0$,

$(s_0s_1s_2s_3s_4)^p s_0 s_1$,

$(s_0s_1s_2s_3s_4)^p s_0 s_1 s_2$,

$(s_0s_1s_2s_3s_4)^ps_0s_1s_2s_3$.

In the another example the block $(s_0s_4s_1s_2s_3)$ can be followed
by itself and the block $s_4s_0s_4$, but we see that the last block
can not be followed by a block beginning with $s_3$ because the
block $s_0s_4s_1s_2s_3s_4s_0s_4s_3$ is the leading term of an
element in the Groebner-Shirshov basis. Then the block
$s_0s_4s_1s_2s_3s_4s_0s_4$ can be followed by a block beginning with
no $s_3$. For example it can be $s_1s_2s_3s_4s_0$. Since this block
can be followed by itself we have the following reduced elements in
the flag subset  $\widetilde{A}_4 / A_4$:

$(s_0s_4s_1s_2s_3)^n$,

$(s_0s_4s_1s_2s_3)^ns_0$,

$(s_0s_4s_1s_2s_3)^ns_0s_4$,

$(s_0s_4s_1s_2s_3)^ns_0s_4s_1$,

$(s_0s_4s_1s_2s_3)^ns_0s_4s_1s_2$,

$(s_0s_4s_1s_2s_3)^n(s_4s_0s_4)(s_1s_2s_3s_4s_0)^p$,

$(s_0s_4s_1s_2s_3)^n(s_4s_0s_4)(s_1s_2s_3s_4s_0)^ps_1$,

$(s_0s_4s_1s_2s_3)^n(s_4s_0s_4)(s_1s_2s_3s_4s_0)^ps_1s_2$,

$(s_0s_4s_1s_2s_3)^n(s_4s_0s_4)(s_1s_2s_3s_4s_0)^ps_1s_2s_3$,

$(s_0s_4s_1s_2s_3)^n(s_4s_0s_4)(s_1s_2s_3s_4s_0)^ps_1s_2s_3s_4$.

\noindent {\bf Acknowledgement:} A part of this work was completed
while the first author were visiting the Mathematical Research
Institute, University Putra Malaysia. Therefore the authors
gratefully acknowledge that this research was partially supported
by the University Putra Malaysia under the e-Science Grant
06-01-04-SF0115 and the Research University Grant Scheme
05-01-09-0720RU.

\end{document}